\documentclass [12pt,a4paper, amsppt, reqno]{amsart}
\evensidemargin
-2mm \oddsidemargin -2mm
 \input amssymb.sty

\DeclareMathOperator{\rank}{rank} 
 
\DeclareMathOperator{\const}{const}

\chardef\bslash=`\\ % p. 424, TeXbook
\newtheorem{Theorem}{Theorem}
\newtheorem*{thm*}{Theorem}
\newtheorem*{t*}{Theorem K}

\newtheorem{cor}{Corollary}

\newtheorem{Conjecture}{Conjecture}

\newtheorem{Lemma}{Lemma }
\newtheorem{property}{Property}
\newtheorem{Proposition}{Proposition}
\newtheorem*{prop*}{Proposition}

\theoremstyle{definition}
\newtheorem{Definition}{Definition}
\newtheorem{Example}{Example}
\newtheorem*{remark*}{Remarks}
\newtheorem*{MainCon*}{Main conclusion}

\newtheorem*{defn*}{Definition}

\newtheorem{Remark}{Remark}

\newcommand{\thmref}[1]{Theorem~\ref{#1}}
\newcommand{\secref}[1]{Section~\ref{#1}}

\newcommand{\lemref}[1]{Lemma~\ref{#1}}

\newcommand{\exampref}[1]{Example~\ref{#1}}

\newcommand{\corlref}[1]{Corollary~\ref{#1}}

\newcommand{\propref}[1]{Proposition~\ref{#1}}
\newcommand{\prref}[1]{Property~\ref{#1}}
\newcommand{\itemref}[1]{(\ref{#1})}

\newcommand{\CL}{\mathcal{L}}

\newcommand{\CO}{\mathcal{O}}

\newcommand{\Frac}{\mbox{Frac}}

\newcommand{\wh}{\widehat}

\renewcommand{\sectionmark}[1]{}

\newcommand{\LND}{\mbox{LND}}
\newcommand{\Pic}{\mbox{Pic}}

\newcommand{\Der}{\mbox{DER}}\newcommand{\DER}{\mbox{DER}}

\newcommand{\1}{^{-1}}

\newcommand{\capl}{\operatornamewithlimits{\bigcap}\limits}
\newcommand{\cupl}{\operatornamewithlimits{\bigcup}\limits}
\newcommand{\suml}{\sum\limits}

\newcommand{\bk}{\bigskip}

\newcommand{\g}{\gamma}
\newcommand{\be}{\beta}

\newcommand{\lm}{\lambda}

\newcommand{\om}{\omega}
\newcommand{\ov}{\overline}
\newcommand{\vp}{\varphi}

\newcommand{\Vp}{\Phi}
\newcommand{\Varphi}{\Phi}
\newcommand{\BC}{\mathbb{C}}
\newcommand{\C}{\mathbb{C}}\newcommand{\BP}{\Bbb{P}}
\newcommand{\BQ}{\Bbb {Q}}

\newcommand{\BN}{\Bbb{N}}
\newcommand{\BZ}{\Bbb{Z}}

\newcommand{\ep}{\epsilon}

\newcommand{\supp} {\operatorname{supp}}

\newcommand{\de}{\partial}

\renewcommand{\a}{\alpha}

\newcommand{\der}{\mbox{DER}}

\begin{document}

\date{}

\title { On $\C$-fibrations over projective curves}

\author {Tatiana Bandman}

\author{Leonid Makar-Limanov}

\thanks{The first author is partially supported by the Ministry
of Absorption (Israel), the Israeli Science Foundation (Israeli
Academy of Sciences, Center of Excellence Program), the Minerva
Foundation (Emmy Noether Research Institute of Mathematics). }

\thanks{The second author is  partially
supported by an NSA grant, Max-Planck Institute of Mathematics,
and Weizmann Institute of Science.}

\address{T. Bandman, Department of
Mathematics, Bar-Ilan University, 52900, Ramat Gan, ISRAEL}
\email{bandman@macs.biu.ac.il}

\address{L. Makar-Limanov, Department of
Mathematics, Wayne State University, Detroit, MI, 48202, USA}
\email{lml@math.wayne.edu}

%  \endsubjclass
% \endkeywords

\maketitle

\begin{abstract}
The goal of this paper is to present a modified version  ($GML$)
of $ML$ invariant which should  take into account  rulings
over a projective base and allow further stratification of 
smooth affine rational surfaces.
We provide a non-trivial 
example where  $GML$ invariant is  computed for a smooth affine rational 
surface admitting no $\BC^+$-actions.  We apply $GML$ invariant to computing $ML$ invariant
of some threefolds.
\end{abstract}

\section{ Introduction }\label{intro }
\baselineskip 20pt \setcounter{equation}{0}
%\numberwithin{equation}{section}

Rational affine surfaces, i.e. affine surfaces birationally
equivalent to a plane is an interesting and rich class of surfaces
worthy of investigation. One of the tools which was used for 
classification of such objects is so called $ML$ invariant    ($ML(S)$)
of a
surface  $S$ which is a characteristic subring of the ring of regular
functions of $S.$ It  consists of the regular functions which are
invariant under all possible $\C^+$-actions on $S$.
Any $\C^+$-action on a surface induces a $\BC$-fibration over an affine curve.
The invariant answers to the question how many fibrations of this kind  the
surface  admits.

Naturally enough the fibrations over projective base  are less  studied  ((\cite{DR},
 \cite  {GM} \cite{Za} \cite{KiKo} \cite{GMMR}).

The goal of this paper is to present a modified version of $ML$
invariant which should  take into account projective rulings
also and
allow further stratification of 
rational surfaces. Of course, it is much easier to introduce an
invariant than to be able to compute it in a particular case. 
We still do not know how to compute $ML$ invariant for a given
surface though some technique is available (see \cite{KML1},\cite{KML2}).

Unfortunately computation of the modified version of the invariant
is even more involved. Nevertheless we present a non-trivial
example where we were able to finish the computation. Hopefully
further techniques will be developed in   due course.

Let us recall the definition of the $ML$ invariant. Let
$R$ be the ring of regular functions of an affine algebraic
variety $V$. Let $LND(R)$ be the set of all locally nilpotent
derivations ($lnd$) of $R$. Then $ML(R) =ML(V)= \bigcap_{\de \in LND(R)}
R^{\de}$ where $R^{\de}$ stands for the kernel of $\de.$

Here is the modified version. Let $F(R) = Frac(R)$ be the field of
fractions of $R$. Take an element $f \in F(R)$ and consider the ring
$R[f] \subset F(R)$, i.e., extension of $R$ by the polynomial
functions of $f$. Call $\de \in \der(F(R))$ a generalized locally
nilpotent derivation ($glnd$)  of $R$ if it is locally  nilpotent on $R[f]$
and $\de(f) =  0$. Define $GML(R) =
\bigcap_{\de \in GLND(R)} F(R)^{\de}$ where $GLND(R)$ is the set of
all generalized locally nilpotent derivations of $R$.

If $R = \CO(S)$, the ring of regular functions on
a surface $S$, we will denote $F(R)$ by $F(S).$
Of course, $F(S)^{\de}$ is the algebraic
closure of $\C(f)$ in $F(S)$ when $\de \in  GLND(R)$. 

Therefore $GML(R)$ is \\
-either $F(R)$  when the only element of $GLND(R)$ is  the zero
derivation, \\
-or a field of rational functions of a curve $C$ when
non-zero $lnd$ are possible on $R[f]$ only for $f \in \C(u)$ where
$u$ is a fixed element of $F(R)$, \\
-or $\C$ when there are at least two
substantially different possible choices of $f$.

If $S$ is rational then $C \cong \BP^1$.

Geometrically speaking if $R = O(S)$ where $S$ is a surface, a non-zero $glnd$ of $R$ which is not
equivalent to an $lnd$ of $R$ corresponds to $\C$-fibration of $S$ over a
projective curve. Therefore $S$ contains a cylinder like subset. By a result of
M. Miyanishi and T. Sugie  (\cite{MiSu} ) it  is equivalent to $\overline{k} = -\infty$
where $\overline{k}$ is the logarithmic
Kodaira dimension of $S$. We can think about $LND(R)$ as a
subset of $GLND(R)$ (just take $f = 1$). So in the case of surfaces the
logarithmic Kodaira dimension of $S$ is $-\infty$ if and only if $GLND(R)$
contains a non-zero derivation.

In  \secref{properties} we give some definitions and
demonstrate the first properties of $GML$\textbf{.}

In
\secref{ex} we compute invariant $GML$  for a ``rigid'' surface:
smooth affine rational surface admitting no $\BC^+$-actions.  In
\secref{sec1} we apply $GML$ invariant to computing $ML$ invariant
of some threefolds.

It appears that $GML$ invariant of a surface $S$ is closely connected  to $ML$
invariant of the line bundles over $S$. Namely, let $\CL=(L, \pi, S)$
be a line bundle over $S$ and $\partial \in LND (\CO(L)).$
Then there exists $\partial '\in GML(S)$ such that  $\partial f=0$
for any $f\in\pi^*(F(S)^{\partial'})$ (see \propref{homo}).
On the other hand for any $\partial\in GML(S)$ there is a line bundle
$\CL=(L, \pi, S)$ and a $lnd$ $\partial' \in LND (\CO(L))$ such
that $\partial' f=0$
for any $f\in\pi^*(F(S)^{\partial})$  (\lemref{gmltoml}).

  This is why the $GML$ invariant is useful for understanding whether
$ML$-invariant of a surface is stable under reasonable geometric
constructions. In our previous work the cylinder over a surface
played the role of a ``reasonable" geometric construction. Here we
are replacing the cylinder by an algebraic line bundle.

%Similarly to the cylinder, the line bundle always admits
%a $\BC$-action along its fibers, but this action is not
%necessarily fixed point free.

It is not always possible to generalize the results known for the
cylinders to this setting. E. g. for ``rigid" surfaces
$$ML(S\times \BC)=\CO(S),$$
but it is not valid for some non-trivial line bundles,
because $GML(S)$ is not trivial.

 In
\secref{sec1} , \corlref{rigid} we describe the line bundles for
which the equality  nevertheless is true.

  Below  we denote by $\BC_{x_1,...,x_n}^n$ the $n$ - dimensional
 complex affine space with
 coordinates $x_1,...,x_n,$ and
 for an irreducible subvariety $C$  of codimension 1 we denote by
  $[C]$
the effective divisor with this support and coefficient 1; $\supp (G)$
and $Cl(G)$ stand for  support and class  of divisor $G$ respectively.
$(C_1,C_2)=([C_1],[C_2])$ is the intersection number of two curves (resp. divisors).
   $\ov{A}$ stands for a closure of $A.$
  {} For a rational function $f$ we denote by $(f),(f)_0,  (f)_{\infty} $
 divisors of $f,$ of its zeros and of its poles respectively.
 If  $\CL=(L, \pi, S)$ is a line bundle over a smooth surface $S$, then $D_L$
 stands for the Weil divisor (since  $S$ is smooth we do not distinguish between Weil and Cartier divisors)
 on $S,$ associated to  $\CL .$
Two $\BC^+$-actions are equivalent if they have the same
generic orbit.

 For a ring $R$ we denote by $\DER(R)$ the set of derivations on $R,$
 by $LND(R)\subset \DER(R)$ the set of locally nilpotent derivations,
 by $F(R)$ the field of fractions of $R.$ For  a derivation  $\de\in \DER (R)$
 we denote by $R^{\de}$ and $F(R)^ {\de}$ the kernel of $\de$ in $R$ and $F(R)$
 respectively.

 The main information on properties of $LND(R)$ may be found in  \cite{KML2}.
 Our Encyclopedia  on affine surfaces with fibrations is the book of M. Miyanishi
 \cite {Mi2}.

\section{Properties of GML}\label{properties}

Let $S$ be a smooth affine complex surface, $R=\mathcal O(S)$ be
the ring of regular functions on $S$, and $F(S)=\Frac(\mathcal O(S))$
stand for the field of fractions of $\mathcal O(S).$

\begin{Definition}\label{glnd} The derivation $\de\in\Der(R)$
is a generalized locally nilpotent
derivation ($glnd$) if there is $f\in F(S),$
 such that $\de\in LND(R[f])$ and $\de(f)=0.$
The set of all $glnd'$s for the ring $R$ is denoted by $GLND(R)$
\end{Definition}

\begin{Definition}\label{eqv} Two elements $\de_1$ and $\de_2 $ in $GLND(R)$
are equivalent, if  $F(R)^{\de_1}=F(R)^{\de_2}.$\end{Definition}

\begin{Definition}\label{GML}
The invariant $GML(R)$ ( or $GML(S)$ if  $R=\mathcal O(S)$)
 is the field $\capl_{\de \in GLND(R)} F^{\de}.$
\end{Definition}

\begin{Definition}\label{GML}
A  smooth affine rational surface  $S$ is rigid, if  $\log$-Kodaira dimension
$\ov k (S)=-\infty$ and $ML(S)=\CO(S).  $
\end{Definition}

The invariant
$GML(S)$ has the following properties:

\begin{property}\label{property1} $\overline k(S)=-\infty$ if and
only if $GML(S)\ne F(S).$
\end{property}

\begin{proof} Indeed, by definition, $GML(S)\ne F(S)$ is
equivalent to the existence of a cylinder-like subset in $S$,
which is equivalent (\cite{MiSu}) to $\overline k(S)=-\infty.$
\end{proof}

\begin{property}\label{property2} If there exists a Zariski open
affine subset $U\subseteq S$ such that $ML(U)=\mathbb C,$ then
$GML(S)=\mathbb C.$
\end{property}

\begin{proof}  Since $ML(U)=\mathbb C,$
the surface $S$ is rational.
Let $\varphi_1:U\to\mathbb C$ and $\varphi_2:
U\to\mathbb C$ be two $\mathbb C$ fibrations on $U.$ Let
$\overline S$ be a closure of $S,$ such that the rational
extensions $\overline\varphi_1:\overline S\to\mathbb P^{1}$ and
$\overline\varphi_2:\overline S\to\mathbb P^{1}$ of $\varphi_1$ and
$\varphi_2$, respectively, are regular. Let $S=\overline S-D,$\
$U=\overline S-(D\cup D'),$  $D=\cupl_{k=1}^m C_k,$\ $D'=\cupl_{j=1}^n B_j, $
where $C_i$ and $B_i$ are reduced irreducible components
of $D$ and $D'$ respectively.  All of them are smooth and rational ( see \cite {Mi2},
ch. III, Lemma 1.4.1).

We denote by $D_1^\infty,D_2^\infty$ such components in $D\cup D'$
that $\overline \varphi_i: D_i^\infty\to\mathbb P^{1}$ is an
isomorphism, $i=1,2.$ 
If $D_i^\infty \subset D'$ for some $i$ then for a
generic $a\in\mathbb P^1$ the intersection
$\overline\varphi_i^{-1}(a)\cap(D\cup D')\in D' - D,$ since $\overline
\varphi_i^{-1}(a)$ has only a single point in $\overline S-U.$ It
follows that a compact curve $\overline\varphi_i^{-1}(a)\subset S,$
thus $S$ is not affine. The contradiction shows that $D_i^\infty\subset D$ for
$i=1,2$ and $\overline\varphi_i\bigm |_{B_j}=\const$ for every
$j=1,\dots,n.$ Thus $\overline\varphi_i\bigm |_{S}:S\to\mathbb P^1$ are
nonequivalent $\mathbb C$-fibrations, and $GML(S)=\mathbb C.$
\end{proof}

\begin{property}\label{property4} (See \cite{GMMR}, \cite{Za}.) For a $\BQ$-
homology plane $S$
$$GML(S)=\Frac(ML(S)).$$
\end{property}

\begin{property}\label{property5} (see \cite {GM}, Th. 4.1) If there exist a $\BC$-fibration
$f:S\to B$ and the curve $B \cong \BC$ ($B \cong \BP^1$), all the
fibers of $f$ are irreducible and there are at least two (resp.
three)  multiple fibers, then $GML(S)=\BC(f).$
\end{property}

 The next Lemma is a simple fact about locally-nilpotent derivations,
  which was proved in another form in \cite{BML1}. We will need it further.

  \begin{Lemma}\label{pine}
Let $R$ be a finitely generated ring and $r \in R$.
Assume that there is a non-zero $lnd$ $\de$ on $R[r^{-1}]$.
 Then there is a non-zero $lnd$ on $R.$
\end{Lemma}
\begin{proof}
Indeed, let $r_1, \dots, r_n$ be a generating set of $R$.
Then $\de(r_i) = p_ir^{-d_i}$ where $p_i \in R$ and $d_i$
is a natural number. It is clear that $\de(r) = 0$ since both
$r$ and $r^{-1}$ are in $R[r^{-1}]$.

Take $m$ which is larger then all $d_i$. Then $\epsilon =  r^m\de$
is also an $lnd$ on $R[r^{-1}]$. Since $\epsilon(r_i) \in R$ for all $i$
the derivation $\epsilon$ is a derivation of $R$. So it is an $lnd$ of $R$.
\end{proof}

\bk

\begin {Remark}\label{rem2}
Same consideration works if there is a non-zero $lnd$ $\de$ on $R(r)$.
Again $\de(r) = 0$ and instead of $r^m$ take a common denominator
of all $\de(r_i)$ which is a polynomial in $r$.
\end{Remark}

\section{Example}\label{ex}

 In this section we compute $ GML(S) $ for a surface $S \subset \BC^7,$
introduced in \cite{BML2} ( example 3) and
 defined by

\begin{equation}\label{e1}
uv=z(z-1)\end{equation}\begin{equation}\label{e2}
v^2z=uw\end{equation}\begin{equation}\label{e3}
z^2(w-1)=xu^2\end{equation}\begin{equation}\label{e4}
u^2(z-1)=tv\end{equation}\begin{equation}\label{e5}
(z-1)^2(t-1)=yv^2\end{equation}
\begin{equation}\label{e6}
u^2v^2=wt\end{equation}\begin{equation}\label{e7}
yz^2=u^2(t-1)\end{equation}\begin{equation}\label{e8}
x(z-1)^2=v^2(w-1)\end{equation}\begin{equation}\label{e9}
v^4x=w^2(w-1)\end{equation}\begin{equation}\label{e10}
u^4y=t^2(t-1)\end{equation}\begin{equation}\label{e11}
v^3=(z-1)w\end{equation}\begin{equation}\label{e12}
u^3=tz \end{equation}\begin{equation}\label{e13}
xy = (w - 1)(t - 1).\end{equation}

Equations \eqref{e6}-\eqref{e13}
 are the consequences of
 the equations \eqref{e1}-\eqref{e5}.

The surface is smooth, because the rank of the Jacobi matrix of
equations \eqref{e1}-\eqref{e13} is maximal everywhere.

The surface $S$  has  the following  properties

\begin{property}\label{A}
\begin{enumerate}
    \item \label{a1} $\ov{\kappa}(S)=-\infty;$
    \item \label{a2} $ R=ML(S)=\CO(S);$
    \item\label{a3} $\pi_1(S)=\BZ/2\BZ$
    \item \label{a4}  $\Pic (S)=\BZ+\BZ/2\BZ;$  
    \item \label{a5} It admits an automorphism $ a: (u,v,z,t,w,x,y)\to (-v,-u,1-z,w,t,y,x);$
    \item  \label{a6} Morphism $b:S\to \BP^1,$ defined as $  b(s)=\frac{z}{u}$ for a point $s\in S,$
    is a $\BC$-fibration. All the fibers of this fibration are isomorphic to $\BC^1.$
    The fibers $B_0=b^{-1}(0)$ and $B_{\infty}=b^{-1}(\infty)$ have multiplicity $2.$
        \item \label{a7}The following relations are valid:
    $$z=ub,  \ v=b(ub-1), \ w=b^3(ub-1)^2, \ $$
    $$ x=b^2(b^3(ub-1)^2-1), \
              t=    \frac{u^2}{b}, \ y=\frac{u^2-b}{b^3}.$$
    \item  \label{a8}The surfaces $S_0=S-B_{\infty}$ and  $S_{\infty}=S-B_0$ are isomorphic to the hypersurface
 $S'=\{\be^3\g=\a^2-\be \}.$ Isomorphisms $\tau_0:S_0\to S'$ and $\tau_{\infty}:S_{\infty}\to S'$ are defined
 respectively by
$\be=b, \a=u,\g=y$ and
    $\be=1/b, \a=v,\g=x.$ Indeed, $R[b]=\BC[u,b,\frac{u^2-b}{b^3}]$ and
    $R[1/b]=\BC[b(ub-1),1/b,b^2(b^3(ub-1)^2-1)].$
    \end{enumerate}
    \end{property}

    \begin{Theorem}\label{GML}
    $GML(S)=\BC(b).$
    \end{Theorem}

    \begin{proof}[Proof of \thmref{GML}]

The proof is rather long but the main idea is as follows: if
$GML(S)\ne\BC(b)$ then there exists a $\BC$-fibration
$\vp\in F(S),$  such that $LND (R[\vp])\ne\{0\},$ where $R=
\CO(S)$ and $\vp$ is algebraically independent with $b$.
    %=
    %\BC[u, ub, b(ub-1), b^3(ub-1)^2, b^2(b^3(ub-1)^2-1),\frac{u^2}{b},\frac{u^2-b}{b^3}]$$
We introduce some weights for the generators $u, v, z, t, w, x, y$
and consider the corresponding graded algebra $\wh {R[\vp]}$
(since $\vp$ is a rational function the weight of $\vp$ will be
also defined). We will show that for these weights $LND
(\hat R)=\{0\}$, where $\hat R$ is a corresponding to $R$ graded algebra.
Then we will get that the leading forms of the numerator and the
denominator of $\vp$ are algebraically dependent, and finally that
$LND (\wh{ R[\vp]})=\{0\}$. This will bring us to a contradiction
because
        %$\om(u)=4, \om(b)=-1, \om(z)=3, \om(v)=2, \om(w)=3,
        %\om(x)=1, \om(t)=9,\om(y)=11.$
(as it was shown in \cite{KML1}, see also \cite{KML2}) $LND
(R[\vp])\ne\{0\}$ implies $LND (\wh{ R[\vp]})\ne\{0\}.$

Let us specify the weights ($\om$) by  $\om(u)  = 4$ and $\om(b) =
-1 + \rho$ where $\rho<<1$ is an irrational number.  Then $\om
(z)=3+ \rho,$  \ $\om (v)=2+ 2\rho,$ \  $\om (w)=3+ 5\rho,$ $\om
(x)=1+ 7\rho,$ $\om (t)=9- \rho,$ $\om (y)=11-3
\rho.$

      %The ring
      %$$\hat R=\BC[u, ub, ub^2, u^2b^5, u^2b^7,u^2b^{-1},u^2b^{-3}]$$

    \begin{Lemma}\label{hat R} $LND (\hat R)=\{0\}$\textbf{.}
        \end{Lemma}
    \begin{proof}[Proof of \lemref{hat R}]
    Let $\de\in LND (\hat R)$ be a non-zero derivation.

    The system

    \begin{equation}\label{f1}
uv=z^2 , \   v^2z=uw , \     z^2w=xu^2, \    u^2z=tv                                               \end{equation}

\begin{equation}\label{f2}
z^2t=yv^2, \ u^2v^2=wt, \ yz^2=u^2t, \  xz^2=v^2w
\end{equation}

\begin{equation}\label{f5}
v^4x=w^3, \ u^4y=t^3, \   v^3=zw, \   u^3=tz, \     xy = wt.
\end{equation}
defines a reduced  (the rank of Jacobian matrix is maximal in Zariski open subset $\{uvz\ne 0\}$)
  and  irreducible surface. The last follows from the fact, that each fiber of
  a rational function $k=\frac{u}{z}=\frac{z}{v}$ is irreducible.

    According to \cite{KML2} (Lemma 6.2)  the system \eqref {f1},\eqref{f2}\eqref{f5}
    defines  $\hat R.$

 So
$${\hat R} = \C[u, z, u^{-1}z^2, u^{-3}z^5, u^{-5}z^7, u^3z^{-1},
u^5z^{-3}].$$  We want to show that this ring does not have a
non-zero locally nilpotent derivation.   With our choice of
weights we will get that the induced non-zero locally nilpotent
derivation $\hat \de$ also belongs to $LND (\hat R)$ since $\hat
R$ is a graded algebra relative to this weights. The weights are
not comesurable, that is why both $\hat \de(u)$ and $\hat \de(z)$
are monomials and $\hat \de$ of any monomial is a monomial.

We present monomials in $u,z$ of $\hat R$ by points of the
two-dimensional integer lattice. The set $A$ of  points
$(1,0),(0,1),(-1,2),(-3,5) ,(-5,7) ,(3,-1),(5,-3)$ which
correspond to the generating set of $\hat R$ is located on a plane
with coordinates $(r,s)$ inside the angle between lines
     $L_1=\{7r+5s=0\}$ and
    $L_1=\{5s+3r=0\}$ containing  the first quadrant. The points $(-5,7),(5,-3)$
    belong to $L_1,$  $L_2$ respectively.
  There is an involution
     $\hat a:(u,z)\to(u^{-1}z^2,z)$ of the ring $\hat R$, which
      changes roles of  lines $L_1$ and $L_2.$

Since $\hat \de$ is locally nilpotent and  non-zero it implies
that there is a monomial $f \in ker(\hat \de) \setminus \C$. This
means that the $ker(\hat \de)$ is generated by a monomial, say
$f$. Now, let us take a monomial $g$ for which $\hat \de(g) \neq
0$ and $\hat \de^2(g) = 0$. It is known (see \cite {KML2}), that
$\hat R \subset \C(f)[g]$. The action of $\hat \de$  is
represented  on the plane $(r,s)$ as the
translation by the vector, corresponding to $g$ toward
  the line passing through the point corresponding to $f.$
  Since both $f, g \in \hat R$ it implies that $f$ must
be represented by a point of the  boundary line, i.e.  $L_1 $ or $ L_2.$
  Because of involution  we may assume that $\hat\de(u^5z^{-3})=0.$
    Let $\deg$ be the degree function induced by $\hat \de$
    (\cite{FLN}). Then $5 \deg u - 3\deg z=0,$ i.e. $\deg u=3n,
    \deg z = 5n$ for some $n\in \BN.$ Since $\hat \de(u)$ is a monomial,
     $n$ should divide $3n - 1$, so $n = 1$.

Now, $\deg(g) = 1$. So one of the monomials in $\hat R$ has degree
1. But $\deg u = 3, \deg z = 5, \deg u^{-1}z^2 = 7, \deg u^{-3}z^5
= 16, \deg u^{-5}z^7 = 20,  \deg u^3z^{-1} = 4, \deg u^5z^{-3} =
0$ and since $g$ is a product of these monomials it cannot have
 degree
equal to 1.
    \end{proof}
    \bk

  The next step is  computation of  the leading form $\hat \vp$ of the function $\vp.$
We need several Lemmas. We will  denote by the same letter the function $u$
 on $S,$ its extension to $\ov S$ and its lift to any blow-up $\tilde S$ of $S.$

    \begin{Lemma}\label{closure}
    The map $b$ can be extended to a morphism $\ov{b}:\ov {S}\to \BP^1$ to the closure $\ov {S}$
    such that   the divisor $D=\ov {S}-S$ has the following graph:

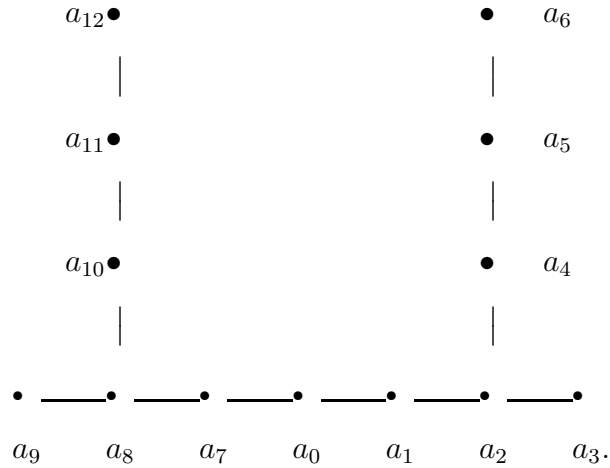
\begin{alignat} {15}
&       &  a_{12}   &   \bullet   &    &    &     &      &      &       &          & \bullet     &  a_6  &  &   {}\notag\\
&       &           &  \bigm|     &    &    &     &      &      &       &          &\bigm |      &       &  &  {}\notag\\
&       &    a_{11} &   \bullet   &    &    &     &      &      &       &          & \bullet     &  a_5  &  &   {}\notag\\
&       &           &  \bigm|     &    &    &     &      &      &       &          &\bigm |      &       &  &  {}\notag\\
&       &    a_{10} &   \bullet   &    &    &     &      &      &       &          & \bullet     &  a_4  &  &   {}\notag\\
&       &           &  \bigm|     &    &    &     &      &      &       &          &\bigm |      &       &  &  {}\notag\\
&{}_\bullet &\underline{\quad\quad}& {}_\bullet   &\underline{\quad\quad}&{}_\bullet&\underline{\quad\quad}&{}_\bullet& \underline{\quad\quad} &{}_\bullet&\underline{\quad\quad}&{}_\bullet &\underline{\quad\quad}&{}_\bullet&
{}\notag\\
&  a_9 &           &a_8      &          & a_7   &    & a_0   &   & a_1   &    &   a_2 &          & a_3 &{}\notag.
\end{alignat}
%%%%%%

%   $$\begin{align}}
%&       & a_8  & a_4    & \\
%&       &\bullet &     \bullet &  \\
%&       &\bigdownline &    \bigdownline &  \\
%   &\underset {a_7} {_\bullet}\underline &\underset {a_6} {_\bullet}&\underline\underset {a_5}
% {_\bullet}\underline\underset {a_0}{_\bullet}\underline
%\underset {a_1} {_\bullet}\underline &\underset {a_2} {_\bullet}&\underline\underset {a_3}
% {_\bullet}\end{align}
% $$
where vertex $a_i, \  0\le i\le 12$ represent a component $A_i$ of divisor $D$.
Moreover, they enjoy

 \begin {property}\label{B}
   \begin{enumerate}
    \item \label{b1} $A_i^2=-2$ for $i>0;$
    \item \label{b2} $A_6\cap B_0\ne\emptyset,$ \ \ $A_{12}\cap B_{\infty}\ne\emptyset;$
    \item \label{b3}$F_0=\ov {b}^{-1}(0)=A_1+A_3+2A_2+2A_4+2A_5+2A_6+2B_0,$ \ \
    $F_i=\ov {b}^{-1}(\infty)=A_7+A_9++2A_8+2A_{10}+2A_{12}+2A_{11}+2B_{\infty};$
    \item \label{b4}$u \bigm|_{\cupl_{2}^{6} A_i}=u \bigm|_{B_0}=u \bigm|_{B_{\infty}}=0,  \ u \bigm|_{A_1}$
    is linear;\\
    $v \bigm|_{\cupl_{8}^{12} A_i}=v \bigm|_{B_0}=v \bigm|_{B_{\infty}}=0,  \ v \bigm|_{A_7}$ is linear.
    \end{enumerate}
\end{property}  \end{Lemma}
\begin{proof}[Proof of \lemref{closure}]

 Due to \prref{A}\itemref{a5} and  \prref{A}\itemref{a8} it is sufficient to analyze the structure of
 the closure of the surface $S_0=S-B_{\infty}$ and to proof only
 \prref{B}, \itemref{b1}-\itemref{b4}.
 The detailed description of the
 graph of the divisor $D_0=\ov {S_0}-S_0$ is given in \cite {MiMa} and \cite{TtD},
 together with the proof of \prref{B}, \itemref{b1}-\itemref{b3}.

  In order to obtain $S_0$ you have to consider the open set $U\cong \BC_{b,u}^2$
 of a Hirzebruch surface and to blow-up several times the point $b=0, u=0$
 of the fiber $B= \{b=0\}.$
That is why \prref{B}\itemref{b4} is   valid: $u=0$ on
all  exceptional components $A_i, 1\leq i\leq
6$ of this process and $u$ is linear along the proper transform
$A_1$ of $B.$ The equality $u \bigm|_{B_{\infty}}=0$ follows from
 equation $\eqref{e1}$ in the definition of the surface.
\end{proof}

Any non-equivalent to $b$ fibration $\vp: S\to  \BP^1, \vp \in F(S)$
 has the following

\begin {property}\label{C}
\begin{enumerate}
\item  \label{c1}every fiber $\Vp_q=\vp^{\1}(q)$ is isomorphic to $\BC$
     (since $\rank \Pic_{\BQ}(S)=1)$ (\cite{Mi2}, Ch.3, 2.4.3.1 p );
    \item    \label{c2}There are precisely two values $q_0, q_1\in\BP^1, $ such that
    the fibers $\Vp_{q_0},$  \ $\Vp_{q_1}$ have multiplicities $2$ ;
     all other fibers are of multiplicities $1$(\cite{Fu}, 4.19. 4.20, 5.9);
        \item  \label{c3}$\vp$ is not a function of $b$ (since they define  non-equivalent
        fibrations);
        \item  \label{c4}  there is     $\de\in LND(R [\vp])$ such that $\de \ne\{0\}$ and $\de (\vp)=0.$
    \end{enumerate}\end{property}

\begin{Lemma}\label{fiber} There is no $p\in\BP^1$ such that $\vp$ is constant
 along the fiber $B_p=b^{\1}(p).$
\end{Lemma}

\begin{proof}[Proof of \lemref{fiber} ] If such $p$ exists, then the affine surface
$S''=S-B_p$ admits two non-equivalent
    $\BC$   fibrations over $\BC, $ i.e. $ML(S'')=\BC$.

    If $p\ne 0,\infty$, then $b\bigm |_{S''}$ has two singular fibers, thus $ML(S'')\ne\BC$
    (\cite {Giz},\cite{Ber}).

    If $p=0$ or $\infty$, then $S''\cong S'$ ( see \prref{A}\itemref{a8}). But $ML(S')=\BC[\be]\ne\BC$
    as well (\cite{MiMa}, Theorem 2.3).

    Thus, both cases are impossible.
\end{proof}

\begin{Lemma}\label{vp}  The extension $\ov\vp$ of rational function $\vp$ to $\ov S$
is not regular and has only one singular point.
\end{Lemma}
\begin{proof}[Proof of \lemref{vp}]

 Assume first that $\ov\vp$ is morphism of $\ov S$ onto $\BP^1.$

 Then for one of the components $A_i$ of divisor $D= \ov S-S$
(see \lemref{closure}) the restriction $\vp \bigm|_{A_{i}}:A_i\to\BP^1$
is an isomorphism.  Due to the existence of the automorphism $a$ of $S$
 ( see \prref{A}\itemref{a5}), we may assume that $0\leq i\leq 6.$

\bk

{\bf Case 1}.  $i=0$. Then  the generic fiber $\ov\Vp_{q}=\ov\vp^{-1}(q)=\ov{\vp^{-1}(q)}
\cong\BP^1$  of $\ov\vp$
intersects $A_0$ transversally. Since  the function $u$
is linear along the generic fiber (\cite{BML2}, Ex. 3)  it has a  simple pole
 along $A_0. $ Since this is the only puncture of $\Vp_{q}$
 and $u\in \CO(S),$  the restriction $u \bigm|_{\ov\Vp_{q}}$has the only simple pole
 at the point $A_0\cup\ov\Vp_{q}.$ But it has zero at every point of
 intersections $\ov\Vp_{q}\cap B_0\ne\emptyset$ and $\ov\Vp_{q}\cap B_{\infty}\ne\emptyset.$
 Hence the number of zeros is at least two. The contradiction shows that $i\ne 0.$
 \bk

{\bf Case 2}. $0 < i\leq 6.$ In this case $\ov\Vp_{q}$ intersects $D$ at a point of $A_i$
only and for a general fiber $u$ is finite at the intersection point (see \prref{B}\itemref{b4}).
 Thus it is finite everywhere in $\ov\Vp_{q},$ hence constant. Since the curve
  $\{u=const\}\not\cong \BC$ in $S$ it is impossible.

\bk

  Thus $\ov\vp$ is regular on $S$ but is not a morphism of $\ov S,$ i.e.
  the singular point of $\ov\vp$ is at the puncture of
   the generic fiber $\Vp_{q}$
   ( or, the same, at the intersection of generic fibers $\ov\Vp_{q}$).
   Since  $\Vp_{q}$ has only one puncture, there is only one singular
    point $s\in\ov S.$
\end{proof}

 Let $\ov {b}(s)=p_0.$
 We may assume that $p_0\ne 0$
 ( due to the involution $a$ we may always change the roles of $0$ and
 $\infty.$)

 Let $\pi: \tilde S\to \ov S$
be a resolution of $\ov\vp,$ i.e. $\pi$
is an isomorphism outside $\pi^{\1}(s).$
Let $\pi^{\1}(s)=\cupl_{0}^{k}E_j,$ where  $E_j$
are exceptional components in $\tilde D=\tilde S -S,$ \
let $\tilde\vp=\ov\vp\circ\pi, $ \ \ $\tilde b=\ov b\circ\pi, $
 let $\tilde A_i$ be  proper transforms of $A_i $ and let
$\tilde\vp\bigm|_{E_0}$ be an isomorphism.Then $\tilde\vp$ has to be constant along
each connected component of $\tilde D-E_0.$

Let $\tilde\Vp_{q}=\tilde\vp^{-1}(q)$ and
$\tilde B_{q}=\tilde b^{-1}(q)$ for a point $q\in \BP^1.$
As above,
$\ov\Vp_{q}=\ov{\vp^{-1}(q)}$ and $\tilde\Vp_{q}=\ov\Vp_{q}$ for the generic $q.$

Consider the connected component $R$ of $\tilde D-E_0 $
 containing the proper transform $\tilde A_0$ of $A_0.$

If $\tilde\vp\bigm|_R=\kappa\in \BP^1,$ then $\tilde\Vp_{\kappa}=\tilde\vp^{-1}(\kappa)=
R\cup C, $ where $C=\ov\Vp_{\kappa}$ is the closure  of $\Vp_{\kappa}$
(this means that $C$ is the only component of
$\tilde\Vp$ that intersects $S$).
\bk

\begin{Lemma}\label{intersection}   $\tilde b(s_1)\ne 0,$ where $s_1=R\cap C.$
\end{Lemma}

\begin{proof}[Proof of \lemref{intersection}]
Assume that $s_1\in\tilde b^{-1}(0).$ We remind that $\pi$ is isomorphism in the neighborhood of
$\tilde b^{-1}(0).$
The point $s_1$ cannot be the intersection point of $\tilde A_0$
 and $\tilde b^{-1}(0),$ since  three components ($C, \tilde A_0, \tilde A_1$) of  the fiber of $\tilde\vp$
cannot intersect at a point (\cite{Mi2}, Ch.3 1.4.1). Thus,  $s_1\in (\cupl_{1}^{6} \tilde A_i )-(\tilde A_0\cap \tilde A_1)$ and $u(s_1)$ is finite. But then $u$ is finite at every point of $C,$ which is impossible.
\end{proof}

\bk
 \begin{Lemma}\label{multiplicity}        The fiber $\Vp_{\kappa}$ has multiplicity
 $2$ in fibration $\vp.$
\end{Lemma}

\begin{proof}[Proof of \lemref{multiplicity}]
Let $\tilde\Vp_{\kappa}=\cupl_{0}^{6}\tilde A_i\cup C\cup R_1,$ where
$R_1$ is the union of other components of $R,$
 and let the corresponding divisor $G$
 of the fiber $\tilde\Vp_{\kappa}=\tilde\vp^{-1}(\kappa)$ be  $G=\suml_{0}^{6} k_i\tilde A_i+\ep C+H,$   ($\supp H=R_1$).

 We want to  prove that $\ep\ne 1.$
We have

$(\tilde A_6, G)=0$ implies $-2k_6+k_5=0,$

$(\tilde A_5, G)=0$ implies $-2k_5+k_6+k_4=0,$

$(\tilde A_4, G)=0$ implies $-2k_4+k_5+k_2=0,$

$(\tilde A_3, G)=0$ implies $-2k_3+k_2=0,$

$(\tilde A_2, G)=0$ implies $-2k_2+k_3+k_4+k_1=0,$

$(\tilde A_1, G)=0$ implies $-2k_1+k_2+k_0=0,$

$(\tilde A_0, G)=0$ implies $k_0 (A_0^2)+k_1+(A_0,\ep C+H)=0.$

It follows that $k_1=\frac{3}{2}k_0$ and $k_0(A_0^2+\frac{3}{2})+(A_0,\ep C+H)=0.$ Since $(A_0,\ep C+H)> 0,$
and $k_0>0$ we have $A_0^2\ne -1.$
Along all components of $G$ except $A_0$ and $C$  the map $\tilde b$ is constant.
If any of them were a $(-1)$ curve, it would be possible to contract it. The new divisor still would
 have normal crossings, because it was obtained by blow-up process from the normal crossing divisor.
Hence we may assume that the only   $(-1)$ curve in $G$ is $C.$  But then
it cannot be of multiplicity $1$ (\cite{Mi2}, Ch.3 1.4.1).
According to \prref{C} \itemref{c2}, multiplicity should be $2.$
\end{proof}

%\bk
%\bk
 \begin{Lemma}\label{hat} $\hat\vp=\hat {y}^k$ for some $k\in \BZ.$
\end{Lemma}

\begin{proof}[Proof of \lemref{hat}]
  By bilinear transformation of $\vp$ we may always
   achieve that  $q_0=\kappa=0,q_1=\infty$ (see  \prref{C} \itemref{c2} ).

    According  to \lemref{closure}
 and  \prref{A} \itemref{a8}
 \begin{equation}\label{s1}
 S-B_{\infty}= S_0=\{b^3y=u^2-b\}  .\end{equation}

Since $\Pic(S_0)=\BZ/2\BZ,$ divisors $2[\Vp_0\cap S_0]\cong 0$
and $2[\Vp_{\infty}\cap S_0]\cong 0.$
This implies that there exist  polynomials $P(u,b,y),$
and $Q(u,b,y),$ such that
$$2[\Vp_0\cap S_0]=(P(u,b,y))_0\cap S_0,$$
$$2[\Vp_{\infty}\cap S_0]=(Q(u,b,y))_0\cap S_0$$
and
$$\vp\bigm|_{S_0}=\frac{P(u,b,y)}{Q(u,b,y)}.$$
On the other hand, $\vp\bigm|_{B_{\infty}}\ne const.$ It follows, that in $S$
 \begin{equation}\label{s2}
\vp=\frac{P(u,b,y)}{Q(u,b,y)}.\end{equation}

We may substitute $u^2$ by $b^3y-b$ into polynomials $P$ and $Q$
and obtain

 \begin{equation}\label{s3}
P(u,b,y)=P_1(y)+uP_2(y,b)+bP_3(y,b),\end{equation}
 \begin{equation}\label{s4}
Q(u,b,y)=Q_1(y)+uQ_2(y,b)+bQ_3(y,b),\end{equation}

Along $B_0$ function $y$ is linear, $u=b=0,$
along the generic fiber $B_p$  we have  $b=p$, $u$
is linear, $y=\frac{u^2-p}{p^3}.$

  For two  generic fibers $B_p =b^{-1}(p)\subset S$ and  $\Vp_q=\vp^{-1}(q)\subset S$
  we denote by $|B_p, \Vp_q|$ the number of
  points in their intersection $B_p\cap\Vp_q$ counted with multiplicities.
  We consider  $B_p$ and  $\Vp_q$
     as reduced curves isomorphic to $\BC.$ Recall that
  $\ov B_p$ and  $\ov\Vp_q$
  are the closures  in $\tilde S$ of $B_p$ and  $\Vp_q$ respectively.

 Let for two generic points $p,q\in \BP^1$
\begin{equation}\label{s6}
    |B_p, \Vp_q|=(\ov B_p,\ov \Vp_q) =N.
\end{equation}

For $q\ne 0,\infty$
\begin{equation}\label{s7}
    |B_0, \Vp_q|=(\ov B_0,\ov \Vp_q) =N/2.
\end{equation}

Let $r$ be the multiplicity of zero of  function $\tilde\vp$
 along $\tilde A_0. $  For  a generic $p$

\begin{equation}\label{s11}
|B_p, \Vp_0|=(\ov B_p,\ov \Vp_0) =(N-r)/2\end{equation}
\begin{equation}\label{s10}
|B_p, \Vp_\infty|=(\ov B_p,\ov \Vp_\infty) =N/2\end{equation}

and

\begin{equation}\label{s8}
 |B_0, \Vp_\infty|=(\ov B_0,\ov \Vp_\infty) =N/4.
\end{equation}

In order to compute
$|B_0, \Vp_0|$ we denote by $B= 2B_0+\suml_{1}^{6} n_i\tilde A_i$ the divisor
of zero fiber $\tilde b ^{-1}(0).$ Due to \lemref{intersection} $\tilde B_0$ intersects
 $\ov \Vp_0$ only inside the surface $S,$ thus for the generic $p$

\begin{equation}\label{s9}
|B_0, \Vp_0|=(\ov B_0,\ov \Vp_0)=\frac{1}{2}(B, \ov \Vp_0)=
 \frac{1}{2}(\ov B_p, \ov \Vp_0)=(N-r)/4.\end{equation}

    Combining \eqref{s11}, \eqref{s9}, \eqref{s3},
    we get
 \begin{equation}\label{h1}   deg P_1(y)=|B_0, \Vp_0|=(N-r)/4\end{equation}
   \begin{equation}\label{h2}   2deg_y P+ deg_u P=|B_z, \Vp_0|=(N-r)/2\end{equation}
   Here $deg_sH$ stand for degree of polynomial $H$ relative to indefinite $s.$

 \bk

    Combining \eqref{s10}, \eqref{s8}, \eqref{s4},
    we get
 \begin{equation}\label{h3}   deg Q_1(y)=|B_0, \Vp_\infty|=N/4\end{equation}
   \begin{equation}\label{h4}   2deg_y Q+ deg_u Q=|B_p, \Vp_\infty|=N/2\end{equation}

   For our weights $\om(u)=1, \om(b)=-1+\rho, \om(y)=11-3\rho, $ it gives

    $\hat P=\hat P_1=\hat{y}^{\frac{N-r}{4}},$ \
            $\hat Q=\hat Q_1=\hat{y}^{\frac{N}{4}},$  \ $\hat\vp=\hat {y}^{-\frac{r}{4}}$
\end{proof}

Now we can prove the Theorem. Were  there a fibration
$\vp$, there would be a non-zero locally nilpotent derivation on
$\wh{R[\vp]}.$ Since the system which defines $\hat R$ and the equation 
$\hat\vp {y}^{\frac{r}{4}} = 1$ again define a reduced
irreducible surface we can conclude that $\wh{R[\vp]} = \hat
R[\hat\vp]$.   But
   %$\hat\vp=\hat {y}^{-\frac{r}{4}}$ and  $\hat{y}=u^2b^{-3}\in \hat R$
    that  is impossible due to \lemref{hat R}
 and \lemref{pine}.   \end{proof}

\begin{Remark}\label {rem1}  The curve $\{y=0\}\subset S$  contains two rational curves.
As it was proven, none of them may be included into a $\BC$ fibration
(compare with \cite{GMMR}, where such curves are called anomalous).\end{Remark}

\begin{Conjecture}\label{con} Let $S$  be a rigid surface which admits a morphism $b:S \to \BP^1$ such that the divisor at infinity built as in \lemref{closure} has the graph which is different from the graph in the Lemma only by the number of vertices in the vertical components of the graph. We would like to conjecture that then \thmref{GML} remains valid.
\end{Conjecture}

\bk

\section{ $ML$ invariant of a line bundle over  a rigid surface}\label{sec1}

 In this section we establish a connection between the   $GML$-invariant
 of a surface and  the $ ML$-invariant of the total space of a line bundle over the surface.
  The computation of ML-invariant is often a very
involved matter even for surfaces and cylinders over surfaces. That is why we  find
it interesting to compute the invariant
 for   threefolds of another type.
 We consider line bundles over rigid surfaces.
 The information on $GML(S)$ appears to be
 very helpful.
\bigskip

 Let us remind  some notions and notations which we use in this section.

 The triple $\CL=(L, \pi, X)$, where $L,X$ are affine varieties and $\pi: L\to X $
 is a morphism defines a line bundle if there is a covering of $X$
 by Zariski open affine subsets $U_\a$ such that $L_\a=\pi^{-1}(U_\a)\cong U_\a\times \BC_{t_\a}$
and in the intersection $L_\a\cap L_\be$ the function
$g_{\a\be}=\frac{t_\a}{t_\be}\in \CO(U_\a\cap U_\be)$ and
does not vanish.

Assume that there are functions $h_\a \in F(U_\a)$ (i.e.rational in $U_\a$)
such that $g_{\a\be}=\frac{h_\a}{h_\be}.$ Let the divisor $D_L$ be such, that $D_L\cap U_\a=(h_\a)\cap U_\a$
 (recall that $(h_\a)$ is the divisor of $h_\a$).
We say that  the  divisor $D_L$ and  its class $[D_L]$ are associated to
the line bundle  $\CL$ and vice versa ( since the surface is smooth, we do not differ between
 Cartier  and Weil divisors). If $h_\a \in \CO(U_\a),$ divisor
$D_L$ is effective.
% Since $H^1(X, \CO)=\{0\}$ for an affine variety $X$ (\cite{Ha}, Ch.3 Th. 3.5),
%every line bundle over affine smooth surface $
%S$ has an effective associated divisor.

The  set of functions $f_\a\in \CO(U_\a)$ such that
$\frac{f_\a} {h_\a}=\frac{f_\be}{ h_\be}$ is a globally defined rational function $f\in
F(X)$ defines the section
of $\CL$ by $t_\a(u)=f_\a(u)$ for a point $u\in U_\a.$
If $D_L$ is effective, it has a section $t_\a(u)=h_\a(u)$ which
vanishes (intersects a zero section) at $D_L.$
The quotient of two sections is a rational function on $X.$
The sheaf $\mathcal F$ of germs of sections of the line bundle is coherent, and the global sections
$\Gamma(\mathcal F)$ form a projective module over $\CO(X)$, which generates  $\mathcal F_x$
at every point $x\in X$  as $\CO_{x,X}$-module (\cite{S}, Theorem 2,\S 45, Prop. 5, 41).

Therefore the ring
 $\CO(L)=\CO(X)[tr_0, tr_1,  ..., tr_K],$
 where  $r_i\in \CO(X)$ and $t$ is rational on $L.$

 This ring naturally admits an  $lnd$ $\de_\pi\in LND(\CO(L))$
 such that $\de_\pi f=0$ if $f\in \CO(X)$ and
 $\de_\pi  t=1.$
 The   $\BC^+$-action $\psi_\pi$ corresponding to $\de_\pi$
 acts along the fibers of $\pi.$
 \bigskip
\begin{Lemma}\label{prop:1b}
Let $X$ be a   smooth affine variety admitting a $\mathbb C^+$-action
$\phi:\BC\times X \to X$. Let $\CL=(L, \pi, X)$ be an algebraic
line bundle  over  $X.$
Then the total space $L$ of $\CL$ admits a $\mathbb C^+$-action
 $\phi':\BC\times L \to L$
 such that the image
$\pi(\Varphi')$ of a general orbit $\Varphi'$ of the action $\phi'$ is a
generic orbit of the action $\phi.$
\end{Lemma}

\begin{proof}[Proof of \lemref{prop:1b}]

 Since the  action  $\phi$ corresponds to a $\de \in lnd(R)$  which is non-zero we can find an element $r \in R = \CO(X)$ such that $\de(r) = p \neq 0$ and $\de(p) = 0$. Put $A = R^{\de}[r]$ and $B = \Frac(R^{\de})[r]=\Frac(R)^{\de}[r]$ (\cite{ML}, Lemma 1 of O. Hadas).

 As we know, $r_i \in B$. Consider the ideal generated by $r_i, i=1, ...,K$ in $B$. Since $B$ is a principal ideal
 domain this ideal is generated by some element $q$.
So we can write $r_i = q \rho_i$ ($\rho_i \in B$) and polynomials $\rho_0, \ldots, \rho_K$ are relatively prime.
Thus we can find $\varsigma_0, \ldots, \varsigma_K \in B$ for which $\sum_i \rho_i\varsigma_i = 1$. Since all elements in $B$ are elements of $A$ divided by elements from $R^{\de}$ it means that we can find elements $\widetilde{\varsigma}_i, i=1, ...K $ in $A$ such  that $\sum_i r_i\widetilde{\varsigma}_i = q \Delta$ where $\Delta \in R^{\de}$. Therefore $tq\Delta \in \CO(L)$. Next, $tr_i = tq \rho_i$. Let $\delta \in R^{\de}$ be a common denominator for the coefficients of all $\rho_i$. We can define now $\widehat{\de}$ by $\widehat{\de}(tq) = \delta$, $\widehat{\de}(r) = \delta\Delta,$ $ \widehat{\de}(r') =0$ for every function $r'\in R^{\de}.$
\end{proof}

\begin{cor}\label{co1}
If $ML(X)=\BC$ and $\CL=(L, \pi, X)$ is an algebraic line bundle 
over
$X $ then $ML(L)=\BC.$
\end{cor}

Our main object of interest is  rigid surfaces.

\begin {Definition}  If the generic orbit  of a
$\BC^+$-action $\vp: \BC_{\lm}\times L\to L$
 on
 the total space of a line bundle $\CL= (L, \pi, S)$ over a smooth affine surface $S$
 is
not contained in a fiber of $\pi$
we  will  call $\vp$
   a skew $\BC^+$-action.\end{Definition}

  \vskip 0.5cm

\begin{Example}\label{skew}
  Define  the projection $\pi:\BC^9\to\BC^7$ by
   \begin{equation}\label {t3}
  \pi(u,v,z,w,x,t,y,s,r)=(u,v,z,w,x,t,y)
  \end{equation}
  and define
  the affine variety $L\subset\BC^9$  by
  equations \eqref {e1}--\eqref{e13} and the following ones:
  \begin{equation}\label {t1}
  su=rz\end{equation}
   \begin{equation}\label {b2}
  s(z-1)=rv.\end{equation}
  Then $\CL=(L,\pi, S)$ is a line bundle over the surface $S$
  defined in \secref{ex} by \eqref{e1}--\eqref{e13}.

Indeed in notations of  \secref{ex} $S=S_0\cup S_{\infty}$ and
  $$\pi^{-1}(S_0)\cong S_0\times \BC^1_{r}, \ s=rb;$$

   $$\pi^{-1}(S_{\infty})\cong S_{\infty}\times \BC^1_{s}, \ r=s/b.$$

   There is   $ \de \in\LND(L)$  that is  defined as

  $$\de s=\de r= 0, \de b=0, \de u=s^mr^{n-m},
  \de z=s^{m+1}r^{n-m-1}, $$
  $$ \de v = s^{m+2}r^{n-m-2},
   \de w=2vs^{m+3}r^{n-m-3},
  \de x=2vs^{m+5}r^{n-m-5},$$
  $$\de t=2us^{m-1}r^{n-m+1},
   \de y=2us^{m-3}r^{n-m+3}.$$

  For any $m\geq 3$ and $n\geq m+5$ this  $lnd $ is well defined
  and provides a skew $\BC^+$-action. Note that this line bundle has a section
  $Z=\{r=u,s=z\}\subset L.$ The   divisor  $D$ of
   intersection $Z$
  with the zero section $Z_0$ is associated to $\CL$ divisor. Let
  $C=\{u=0,\  b\ne0, \ b\ne\infty\}$ and let $F$ be a fiber $b=const\ne 0, \infty.$ Then
    $D=C+B_0$ and  since   $(u)_0=C+B_0+2B_\infty\sim 0,$ we have $D\sim -2B_\infty\sim -F.$
  \end{Example}

 \vskip 0.3 cm

  Similar example may be constructed over any  rigid surface $S$.

\begin{Lemma}\label{gmltoml}
   Let $S$ be a rigid surface and $\partial \in GML(S).$ There exists a line bundle
   $(\CL, \pi, S)$
and  $\partial' \in LND (\CO(L))$ such that
$\partial' f=0$
for any $f\in\pi^*(F(S)^{\partial})$ (as we mentioned in Introduction).\end{Lemma}

\begin{proof}[Proof of \lemref{gmltoml}]

Consider  a $\BC$-fibration $f: S\to \BP^1$ on $S$ induced by
$\partial $
 and a non-singular fiber $F=f^{-1}(\infty).$
Consider the line bundle $(\CL, \pi, S)$ associated to the divisor $-mF.$
Let $U_1=S-F$ and $U_2=S-F',$ where $F'=\{f=0\}$ is another non-singular fiber.
( We may always assume that fibers $F$ and $F'$ are nonsingular). Then $L=L_1\cup L_2$
where  $ L_1=\pi^{-1}(U_1)\cong U_1\times \BC_{t_1}$ and
$ L_2=\pi^{-1}(U_2)\cong U_2\times\BC_{t_2}$ and  \ $ t_2=f^m t_1.$
The function $\tau=t_1=f^{-m}t_2\in \CO(L)$: it has zero of order $m$ along $F$
because $f$ has a simple pole there.
The divisor  $(\tau)=Z_0+m\pi^{*}F.$
Thus, $ \CO(L)=\CO(S)[\tau, \tau\om_1^{*}, ...\tau\om_n^{*}],$
where $\om_i$ are rational functions on $S$, such that $(\om_i)\ge -mF.$
  Since $f(U_1)$ is an affine curve, there exists an $lnd$ $\de_1\in LND(\CO(U_1))$  such that
$\de_1 f=0.$ Let $N$ be bigger than the order of poles of $\de_1 \om_i$ along $F$ for all $i=1,...,n.$
One can  define an $lnd$ $\de '\in \CO(L)$ by $\de ' \tau=\de' f=0,
\de ' u=\tau^N\de_1 u $ for $u\in \CO(S).$\end{proof}

\vskip 0.3cm

Take now any   morphism $f: S\to\BP^1$  of a rigid surface $S$
  onto $\BP^1$ such that
   the general fiber of $f$
    is isomorphic to $\BC.$
  Picard group of $S$ is generated by
 divisor $[F]$ of the generic  fiber $F$ and the divisors  $[E_{i,j}]$ of the irreducible components
$E_{i,j}$ of the singular fibers $F_i$, \
 $[F_i]=\suml_{1}^{n_i} \a_{i,j}[E_{i,j}], i=1,..,n   $
with relations reflecting that all the fibers are equivalent.

The group $\Pic(S)\otimes \BQ\cong \BQ^{\oplus N},$
 where $N=(\sum n_i) -n+1$,
 and is generated by  $[F]$ and $[E_{i,j}], j>1.$
 (\cite{Mi2}, Ch.3,  Lemma 2.4.3.1).

 Any element $l\in \Pic(S)$ may be represented uniquely   as
 $$l=m[F]+\suml_{1}^{n}\suml_{1}^{n_i} m_{i,j}[E_{i,j}],$$
 where
 \begin{enumerate}
% \item  For any $i, 1\le i\le n$ at least one of $m_{i,j}=0;$
 \item   $m_{i,j}<\a_{i,j}$ for any $i, 1\le i\le n$ and any $j, 1\le j\le n_i;$
 \item  $m_{i,j}\ge 0$  for at least one of $j, 1\le j\le n_i $
for any $i, 1\le i\le n;$
 \end{enumerate}
 \bk
 \begin{Definition}\label{standard}
 We will call the representation with properties (1)-(2) standard for fibration $f.$
 We will call the element  $l\in \Pic(S)$
 positive relative to fibration  $f$, if  in the standard representation $m \ge 0.$
 \end{Definition}

The crucial fact for  \lemref{gmltoml}  and \exampref{skew}
is that the line bundles  are  associated with  the  non-positive
(relative to a given fibration) element of the Picard group.
 The following example presents the line bundle associated to
a positive divisor.

\bk

\begin {Example}\label{noskew}
  Define  the  same  projection $\pi:\BC^9\to\BC^7$ by
   \begin{equation}\label {n1}
  \pi(u,v,z,w,x,t,y,s,r)=(u,v,z,w,x,t,y)
  \end{equation}
  and
  the affine variety $L\subset\BC^9$  by
  equations \eqref {e1}--\eqref{e13}  and the following ones:
  \begin{equation}\label {n2}
  su=rv\end{equation}
   \begin{equation}\label {n3}
  st=ru(z-1),\end{equation}
   \begin{equation}\label {n4}
  svz=rw.\end{equation}

  Then $\CL=(L,\pi, S)$ is a line bundle over the surface $S,$
  defined in \secref{ex} by \eqref{e1}--\eqref{e13}.

  In notations of  \secref{ex}
  $$\pi^{-1}(S_0-C)\cong (S_0-C)\times \BC^1_{r}, \ s=rv/u    ;$$

   $$\pi^{-1}(S_{\infty}-C_1)\cong (S_{\infty}-C_1)\times \BC^1_{s}, \ r=su/v.$$

Here  $C_1=\{v=0, \ b\ne0, \ b\ne\infty\}$  does not intersect $C=\{u=0,\  b\ne0, \ b\ne\infty\}.$

The associated to $\CL$ divisor being the intersection divisor of the
section $Z_1=\{r=u,s=v\}\subset L$ and the zero section $Z_0$
 is $B_0+B_\infty.$ Therefore $\CL$ is associated to a   positive (relative to the fibration) divisor.  We will show 
that there is no skew
 actions on $L.$
\end{Example}

\begin{Proposition}\label{homo} Let $ \CL=(L, \pi, S)$ be an
algebraic line bundle
over a rigid surface $S.$
Assume that $L$ admits a skew  $\BC^+$-action   $\a:\BC\times L\to L$.
Then there exists a skew $\BC^+$-action  $\be:\BC\times L\to L,$ $\de '\in GML(S)$
and  an induced  by $\de '$ morphism $g:S\to \BP^1$ of $S$
such that
\begin{enumerate}
    \item\label{q1}  the generic fiber of $g$ is $\BC;$
 \item  \label{q2} $g(\pi(O))$ is a point for a general orbit $O$
of $\be;$
 \item  \label{q3} there is no non-zero section $Z$ of $\CL$ over an open subset $U\subset S, $
 such that
 \begin{enumerate}
 \item\label{q31} $g(U)=\BP^1;$
 \item\label{q32} the components of $g^{-1}(p)\cap U$ are isomorphic to $\BC$ for each $p\in\BP^1;$
 \item\label{q33} $g(Z\cap Z_0 )$ is a finite set in $\BP^1.$
  \end{enumerate}
  \end{enumerate}
\end{Proposition}

\begin{proof}[Proof of \propref{homo}]

\bk

\begin{Lemma}\label{lml1}
Let $R$ be an affine ring and $Q= R[t,t^{r_1}\omega_1, \dots,
t^{r_k}\omega_k]$ where $t$ is a variable and $\omega_i \in
\Frac(R)$ . Let $\de \in LND(Q)$ which is not identically zero on
$R$. Then there exists a locally nilpotent derivation on
$Q$ which is $t$-homogeneous and is not identically zero on $R$.
\end{Lemma}
\begin{proof}[Proof of \lemref{lml1}]
Let us introduce a weight function on $Q$ by $w(t)
= 1$, $w(r) = 0$ for $r \in R^*$, and $w(0) = -\infty$.
Consider a (non-zero) locally nilpotent derivation $\overline{\de}$
which corresponds to
this weight function
(\cite{KML2}).
Clearly $\overline{\de} \in LND(Q)$ since
$Q$ is a graded algebra relative to the introduced weight
function. Then $\overline{\de}(t) = t^{k+1}\ep(t)$,
$\overline{\de}(r) = t^k\ep(r)$ where $\ep \in \Der(Q)$ such that
$\ep(t), \ep (r) \in \Frac(R)$ if $r \in R$. Since our goal
is to
produce a locally nilpotent derivation on $R$ we may assume that $k > 0$
  (otherwise
$\de$ can be restricted on $R$). It remains to show that
$\overline{\de}$ is not identically zero on $R$. So assume that
$\overline{\de}$ is identically zero on $R$. Then
$\overline{\de}(t) = t^{k+1}\ep(t)$ implies that
$\overline{\de}(t) = 0$, so $\overline{\de}$ would be identically
zero contrary to the facts. Indeed, if $\deg$ is the degree
function induced on $Q$ by $\overline{\de}$ we have $\deg(t) - 1 =
(k+1)\deg(t) + \deg(\ep(t))$. But since we assumed that
$\overline{\de}$ is identically zero on $R$ we have $\deg(\ep(t))
= 0$ if $\ep(t) \neq 0$. (If $\ep(t) = 0$ then $\deg(\ep(t)) =
-\infty$.) So if $\ep(t) \neq 0$ then $\deg(t) - 1 =
(k+1)\deg(t)$. Since $k > 0$ we see that then $\deg(t) < 0$ which
is impossible. So the lemma is proved.
\end{proof}

\bk

\begin{cor}\label{colml}
 If $\dim(R) > 1$ then
$\Frac(Q)^{\overline{\de}}$ contains a non-constant rational
function from $\Frac(R)$.
\end{cor}

\begin{proof}[Proof of \corlref{colml}] Since $\overline{\de}$ is $t$-homogeneous
the ring
of $\overline{\de}$-constants is generated by $t$-homogeneous
elements. Since $\dim(Q) > 2$ there are two algebraically
independent homogeneous $\overline{\de}$-constants, say $f_1 =
t^m\omega_1$ and $f_2 = t^n\omega_2$. Then
$f_1^n f_2^{-m}\in \Frac(R).$
\end{proof}

\vskip 0.3 cm

We apply the \corlref{colml}
  assuming $R=\CO(S),$ and $Q=\CO(L)$
  and $t\in\CO(L)$ is any regular function on $L$  that is linear along the generic
   fiber and vanishing at zero section.
   Let $\be$ be the $\BC^+$-action defined by locally nilpotent
   derivation
  $\overline{\de}.$
  By construction, all the points
  of zero section $Z_0\subset L$
   are fixed by $\be$ and there exists  $\be$-invariant function
   $f=\pi^*g\in \Frac(\CO(L))$
 with  $g\in \Frac(\CO(S)).$ Using Stein factorization we may assume that the generic fiber of
$g^{-1}(p), p\in \BP^1$ is connected (and irreducible).
 \begin{Lemma}\label{m}
     $g: S\to\BP^1$
   is morphism.
   \end{Lemma}
 \begin{proof}[Proof of \lemref{m}]
  We will identify $S$ with the zero section $Z_0,$  i.e.
  $S\subset L.$ By construction it is $\be$-invariant.
  The function $f$ is the composition of rational maps: $L\overset{\pi}\to S \overset{g} \to \BP^1.$
     Let $p$ be a point in $\BP^1.$
    Let  $C_p=g^{-1}(p)\subset S$ and let   $T_p=\pi^{-1}(C_p)=f^{-1}(p)$.
     Since $f$ is $\be$-invariant,  $T_p$ is $\be$-invariant  as well, thus consists  of $\be$-orbits.
     Since $\be$ is a skew action, these orbits are not mapped to a point by $\pi.$
     Hence, $C_p =\pi(T_p)=T_p\cap Z_0\cong \BC.$ By construction  $T_p$  is the restriction of
     our line bundle $\CL$ over $C_p ,$ thus $T_p\cong \BC^2.$

      If $g$ were not a  morphism there would be  a point $s\in S$
    contained in every fiber  $C_p=\{g=p\}.$ Then for  every $p$
    the set $T_p$ would contain two $\be$-invariant
     intersecting curves:   $C_p$
   and  $A_s=\pi^{-1}(s).$   But then all the points of $T_p$ for all $p$
   would be fixed by $\be.$ The contradiction shows that such point $s$ does not exist
  and $g$ is morphism.
   \end{proof}

   Items \itemref{q1},\itemref{q2} of \propref{homo} are proved in \lemref{m}.
  Assume now that there exists a section $Z$ as in item \itemref{q3}.

   Items \itemref{q31},\itemref{q32},\itemref{q33} imply that $Z\cong U$ admits a $\BC$-fibration
   over $\BP^1$ such that $Z\cap Z_0$ is the union of finite set of fibers of this fibration.
   We want to show that $Z$ is $\be$-invariant and this fibration should be induced by the
   restriction of $\be$ on $Z.$ It would lead to a contradiction, because a $\BC^+$-action has
    an affine base (\cite{MiMa1}, Lemma 1.1).

In notations of \lemref{m}
  item \itemref{q33} means that for a generic $p\in \BP^1$ the curve
   $ B_p=Z\cap T_p$ does not intersect $Z_0, $ in particular, the curves $B_p\subset T_p$ and $C_p=Z_0\cap T_p\subset T_p$ do not intersect.

   Since  $C_p=\pi(T_p),$ and $Z$ is a section, $B_p=Z\cap T_p$  is a section of the bundle
   over $C_p$ and
    $\pi \bigm|_{B_p}:B_p\to C_p$ is an isomorphism. Hence $B_p\cong \BC.$ Thus, in $\be$-invariant set $T_p\cong\BC^2$
    we have two rational
    disjoint curves. $C_p$ is a $\be$-orbit in $T_p,$ therefore the same should be true for
    $B_p.$
    Therefore,  $Z$ is
   $\be$-invariant, and the base of the restriction of the induced fibration should be affine. This contradicts
to \itemref{q31}.
   \end{proof}
 \vskip 0.3 cm

 \begin{cor}\label{rigid} Let $S$ be a
 rigid surface, let $GML(S)=\BC(f)$ and let $f:S\to \BP^1$ be the corresponding fibration.
  Let  $\CL=(L,\pi, S)$
 be a line bundle over $S$. Then $ML(L)=\CO(S)$
 if $\CL$ is associated  to positive (relative to fibration $f$) element $l$ of
 $\Pic(S).$
 \end{cor}

 \begin{proof}[Proof of \corlref{rigid}]
Let $\psi:\BC\times L\to\BC$ be a skew action
 on $L$. According to \propref{homo} it gives rise to  a $\BC$-fibration
 $g:S\to\BP^1.$ Since $GML(S)=\BC(f),$ the fibrations $g$ and $f$ have to be equivalent.
 Let the associated to $\CL$
  element $l\in \Pic(S)$ have the standard representation
 $$l=m[F]+\suml_{1}^{n}\suml_{1}^{n_i} m_{i,j}[E_{i,j}]$$
and let $l_{+}$  be a sum of  summands with non-negative coefficient and
$l_{-}$  be a sum of  summands with negative coefficient. Let $D_{+}$ and $D_{-}$ be
 the union of components appearing in  $l_{+}$ and $l_{-}$ respectively.

  Over $U=S-D_{-}\subset S ,$ the line bundle $\CL$  is associated to  the effective divisor,
  hence has a section $Z_U$ such that intersection $Z_0\cap Z_U\subset D_{+}.$
  Since $\supp D_{+}$ contains at least one component of every fiber of $g, $
   $U$ enjoys all the properties of item \itemref{q3} of
   \propref{homo} which is impossible
if $\psi$ is a skew $\BC^+$-action.
\end{proof}

  \corlref{rigid} provides the
 situation when    similar to
  the case of trivial line bundle,
   the isomorphism  $ML(S)\cong\CO(S)$ implies
     $ML(L)\cong ML(S)\cong \CO(S)$ (\cite{BML2}).
  The following questions remain open.

  {\bf Questions}.

  {\bf 1.}  Let $S$ be a rigid surface, $GML(S)=\BC(f)$
 and let $f:S\to \BP^1$ be the corresponding fibration.
  Let  $\CL=(L,\pi, S)$
 be a line bundle over $S$. Is it possible that  $ML(L)=\CO(S)$
 if $\CL$ is associated  to a non-positive relative to fibration $f$ element $l$ of
 $\Pic(S)?$

{\bf 2.} Assume that  $S$ is rigid and $GML(S)=\BC.$ When $ML(L)=\CO(S)$?

{\it Acknowledgments.} We are grateful to Sh. Kaliman for his
valuable advice and infinite help.
 This
project was started while both authors were visiting Max Planck
Institute of Mathematics and  while the first author
was visiting Wayne State University and University of California, San Diego.

\baselineskip 15pt

\end{document}